\documentclass[12pt,a4paper]{article}
\usepackage{amssymb}
\usepackage{amsfonts}
\usepackage{graphicx}
\usepackage{amsmath}
\usepackage{amsthm}
\usepackage{epsfig}
\usepackage{hyperref}
\usepackage[top=3cm, bottom=2cm, left=2.3cm, right=2.3cm]{geometry}

\newcommand{\lab}[1]{\label{#1}}                

\newtheorem{thm}{Theorem}
\newtheorem{cor}[thm]{Corollary}

\newtheorem{lem}[thm]{Lemma}
\newtheorem{defn}[thm]{Definition}
\newtheorem{cl}[thm]{Claim}
\newtheorem{prop}[thm]{Proposition}
\newtheorem{obs}[thm]{Observation}
\newtheorem*{rem*}{Remark}
\def\qed{\ifhmode\unskip\nobreak\fi\hfill\ifmmode\square\else$\square$\fi}

\def\ss{\smallskip}
\def\no{\noindent}
\def\inst#1{$^{#1}$}
\newcommand{\PW}{PW_2}

\long\def\onefigure#1#2#3#4{
\begin{figure}[ht]
\begin{center}
\scalebox{#2}
{
\includegraphics{#1}
}
\end{center}
\caption{\label{#4} #3}
\end{figure}
} 

\date{}
\title{On the Geometric Ramsey Number of Outerplanar Graphs
\thanks{An extended abstract of this paper appeared in the proceedings of the EuroComb 2013 conference. 
This research was started at the 2nd Eml\'{e}kt\'{a}bla Workshop held in Gy\"{o}ngy\"{o}starj\'{a}n, January 24--27, 2011.
Research was supported by the project CE-ITI (GA\v{C}R P202/12/G061)
of the Czech Science Foundation and
by the grant SVV-2013-267313 (Discrete Methods and Algorithms).
Josef Cibulka and Pavel Valtr were also supported by the project no.\ 52410 of the
Grant Agency of Charles University. Pu Gao was supported by the
Humboldt Foundation and is currently affiliated with University of Toronto.
Marek Kr\v{c}\'{a}l was supported by the ERC Advanced Grant No.~267165.}
}

\author{Josef Cibulka\inst{1}, Pu Gao\inst{2}, Marek Kr\v{c}\'al\inst{1}, Tom\'a\v{s} Valla\inst{3}, Pavel Valtr\inst{1}
}

\begin{document}

\maketitle

\begin{center}
{\footnotesize
\inst{1} 
Department of Applied Mathematics,
Charles University, Faculty of Mathematics and Physics, \\
Malostransk\'e n\'am.~25, 118~00 Praha 1, Czech Republic.\\
\texttt{\{cibulka,krcal\}@kam.mff.cuni.cz}
\\\ \\
\inst{2}
Max-Planck-Institut f\"{u}r Informatik,\\
Saarbr\"{u}cken, Saarland, Germany. \\
\texttt{janegao@mpi-inf.mpg.de}
\\\ \\
\inst{3}
Czech Technical University, Faculty of Information Technology, \\
Prague, Czech Republic. \\
\texttt{tomas.valla@fit.cvut.cz}
}
\end{center}  

\begin{abstract}
We prove polynomial upper bounds of geometric Ramsey numbers of
pathwidth-$2$ outerplanar triangulations in both convex and general cases.
We also prove that the geometric Ramsey numbers of the ladder graph on $2n$ vertices
are bounded by $O(n^{3})$ and $O(n^{10})$, in the convex and general case, respectively.
We then apply similar methods to prove an $n^{O(\log(n))}$ upper bound 
on the Ramsey number of a path with $n$ ordered vertices.
\end{abstract}

{\bf MSC codes:} 52C35, 05C55, 05C10

{\bf Keywords:} Geometric Ramsey theory, Outerplanar graph, Ordered Ramsey theory, Pathwidth

\section{Introduction and basic definitions}

A finite set $P \subset \mathbb{R}^2$ of points is in a \emph{general position}
if no three points of $P$ are collinear. The \emph{complete geometric graph on
$P$}, denoted by $K_{P}$, is the complete graph with vertex set $P$, whose
edges are drawn as the straight-line segments between pairs of points of $P$.

The set of points $P$ is in \emph{convex position} if $P$ is the set of vertices of
a convex polygon. If $P$ is in convex position, we say that $K_{P}$ is a
\emph{convex complete geometric graph}.

K\'arolyi, Pach and T\'oth \cite{KPT} introduced the concept of Ramsey numbers
for geometric graphs as follows. Given a graph $G$, the \emph{geometric
Ramsey number} of $G$, denoted by $R_g(G)$, is the smallest integer $n$ such that every complete
geometric graph $K_P$ on $n$ vertices with edges arbitrarily coloured by two colours
contains a monochromatic non-crossing copy of $G$.
The \emph{convex geometric Ramsey number} of $G$, $R_c(G)$, is defined the same way
except that $K_P$ is restricted to the convex complete geometric graph.
A graph $G$ is said to be {\em outerplanar} if $G$ can be drawn in the plane without any
edge crossing and with all vertices of $G$ incident to the unbounded face.
Apparently, the numbers $R_g(G)$ and $R_c(G)$ are finite only if $G$ is outerplanar:
consider a planar but not outerplanar graph $G$, then it is easy to see that
one cannot find a non-crossing monochromatic copy of $G$ in a convex complete graph.
Also, it follows immediately from the definitions that $R_c(G) \leq R_g(G)$
for every outerplanar graph $G$.

The Ramsey numbers of outerplanar graphs, as well as of all planar graphs,
are bounded by a function linear in the number of vertices 
by a result of Chen and Schelp~\cite{ChenSchelp}.
In contrast, the only known general upper bound on the geometric Ramsey numbers
of outerplanar graphs is exponential in the number of vertices.
This bound follows from the exponential upper bound on the Ramsey numbers for cliques
since a monochromatic clique on $n$ points implies a monochromatic non-crossing occurrence of every
outerplanar graph on $n$ vertices by the result of Gritzmann et al.~\cite{GritzmannEtAl}
(see Lemma~\ref{l:gritzmannEtAl}).

The geometric Ramsey numbers of some outerplanar graphs are known to be both larger than linear
and smaller than exponential,
and it remains open whether there is a general polynomial bound for all outerplanar graphs.
By a simple constructive proof, it is easy to see that for every $n\ge 3$, the cycle graph $C_n$ 
on $n$ vertices satisfies $R_c(C_n)\ge (n-1)^2+1$.
Balko and Kr\'{a}l~\cite{BK} constructed colourings that improve this bound to 
$R_c(G)\ge 2(n-2)(n-1)+2$.
This bound is tight both in the convex and general geometric setting by an earlier 
result $R_g(C_n)\le 2(n-2)(n-1)+2$ of K\'{a}rolyi, Pach, T\'{o}th and Valtr~\cite{KPTV},
This shows that one cannot have geometric Ramsey numbers for general outerplanar graphs
asymptotically smaller than $\Omega(n^2)$.
K\'{a}rolyi et al.~\cite{KPTV} found the exact value $R_c(P_n)=2n-3$
and the upper bound $R_g(P_n) \in O(n^{3/2})$, 
where $P_n$ is the path on $n > 2$ vertices.
The bounds $2n-3 \le R_g(P_n)\le O(n^{3/2})$ remain the best known bounds 
on the geometric Ramsey number of paths.
Further results and open problems on the geometric Ramsey numbers can be found in the survey of 
K\'{a}rolyi~\cite{Ksurvey}.

The ladder graphs are defined as follows.

\begin{defn}\lab{defLadder} For any integer $n\ge 1$, the \emph{ladder graph} on $2n$ vertices, denoted
by $L_{2n}$, is the graph composed of two paths $(u_i)_{i=1}^{n}$ and $(v_i)_{i=1}^{n}$,
together with the set of edges $\{u_i v_i: i \in [n]\}$. See an example in Fig.~\ref{fig:exladder}.
\end{defn}

\onefigure{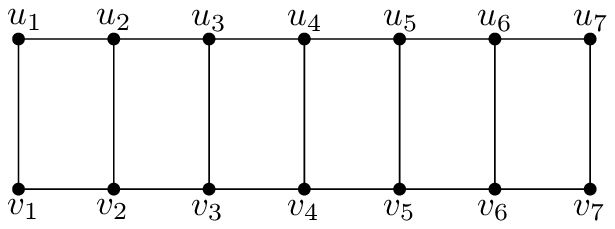}{1}{The ladder graph $L_{14}$.}{fig:exladder}

In this paper,
we contribute to this subject by showing polynomial upper bounds on the geometric Ramsey numbers
of the ladder graphs, and their generalisation.
In Section~\ref{s:ladder}, we show that the geometric Ramsey numbers of the ladder graph on $2n$ vertices
are bounded by  $32 n^3$ and $O(n^{10})$ in the convex and general case, respectively.
In Section~\ref{s:generalisation}, we generalise the polynomial upper bounds to the
class of all subgraphs of pathwidth-2 outerplanar triangulations, see Definition~\ref{def:PW2}.
These bounds are $20 n^7$ and $O(n^{22})$ in the convex and general case, respectively

In Section~\ref{s:ladder:ordered}, we consider the closely related area of the ordered Ramsey theory.
The ordered Ramsey theory recently gained a lot of attention~\cite{FPSS,MS,MSW,BK}, 
mainly in the more general hypergraph setting.
An \emph{ordered graph} $G$ is a graph with a total order $\prec$ on the vertices of $G$.
We say that an ordered graph $G$ is a subgraph of an ordered graph $H$ if the vertices of $G$ can be injectively mapped 
to the vertices of $H$ while preserving both the ordering and the edges of $G$.
The \emph{ordered Ramsey number} $R_o(F,G)$ of ordered graphs $F$ and $G$ is the smallest number $N$ 
such that every $2$-colouring of the edges of the ordered complete graph $K_N$ on $N$ vertices either contains a blue copy of $F$
or a red copy of $G$.

The proof of the upper bound on the convex geometric Ramsey number of the ladder graph in Section~\ref{s:ladder:conv} 
can be extended to show that the ordered Ramsey number $R_o(L_{2n}, L_{2n})$ of the ladder graph $L_{2n}$ with specifically ordered 
vertices is at most $32 n^3$.
The ideas of the proof are applied in Section~\ref{s:ladder:ordered} to give an $n^{O(\log(m))}$ 
upper bound on the ordered Ramsey number $R_o(K_n, P_m)$, 
where $K_n$ is the ordered complete graph on $n$ vertices 
and $P_m$ is an arbitrarily ordered path on $m$ vertices.

We note here that all colourings in this paper,
unless specified, refer to edge colourings.
As a convention, in any $2$-colouring, we assume that the colours used are blue and red.

When $c$ is a colour, we say that $v$ is a \emph{$c$-neighbour} of $u$ 
if the edge $\{u,v\}$ has colour $c$.
Let $N_c(v)$ be the set of $c$-neighbours of a vertex $v$.
We abbreviate the set $\{1,2,\dots, k\}$ with $[k]$ and $\{l,l+1,\dots, k\}$ with $[l,k]$.
We write $(x_i)_{i=1}^{k}$ for the sequence $x_1, x_2, \dots, x_k$.
The sequence of vertices $(v_i)_{i=1}^{\ell+1}$ is a \emph{path of colour $c$} and length $\ell$
if every pair $\{v_i v_{i+1}\}$, $i \in [\ell]$ is an edge and has colour $c$.
A sequence $(A_i)_{i=1}^{k}$ is said to be a \emph{partition} of $A$ if $A_i$ are
pairwise disjoint and $\cup_{i=1}^k A_i = A$.

\section{Ladder graphs}\lab{s:ladder}

In Subsections~\ref{s:ladder:conv}~and~\ref{s:ladder:general}, we prove upper bounds on the convex and geometric 
Ramsey numbers $R_c(L_{2n})$ and $R_g(L_{2n})$ of ladder graphs $L_{2n}$.
Both proofs use the following lemma due to Gritzmann et al.~\cite{GritzmannEtAl}.

\begin{lem}[Gritzmann et al.~1991~\cite{GritzmannEtAl}]\lab{l:gritzmannEtAl}
Let $G$ be an outerplanar graph on $n$ vertices and let $P$ be a set of $n$ points
in general position. Then $K_{P}$ contains a non-crossing copy of $G$.
\end{lem}

In Subsection~\ref{s:ladder:ordered}, a small change to the proof of the upper bound on the convex Ramsey number is shown to 
give an upper bound on the ordered Ramsey number of paths.

\subsection{Convex position} \lab{s:ladder:conv}

\begin{thm}\lab{t:Ladder} For every $n\ge 1$,
$R_c(L_{2n})\le 32n^3$.
\end{thm}

In this section, let $C$ denote a set of $32 n^3$ points in convex position.
That is, $C$ is the set of vertices of some convex polygon.
We label the vertices $v_1, v_2, \dots, v_{|C|}$ in the clockwise order
starting at an arbitrarily chosen vertex $v_1$.
We write $v_i\prec v_j$ if and only if $i < j$.
Let $A, B \subset C$.
We say that $A$ precedes $B$ and write $A\prec B$
if and only if for every $u\in A$ and every $v\in B$, $u \prec v$.
Notice that if $A \prec B$, then the sets $A$ and $B$ can be separated by a line.

For a pair of disjoint vertex sets $(L,R)$, $L \subset C$, $R \subset C$, the complete bipartite graph
on $(L,R)$, denoted by $K_{L,R}$, is the set of edges $\{u,v\}$, where $u \in L$ and $v \in R$.
A complete bipartite graph $K_{L,R}$ is said to be \emph{well-split} if
$L \prec R$ or $R\prec L$.
A well-split $K_{m,n}$ is a well-split $K_{L,R}$, for some $L$ and $R$ such that  $|L| = m$, $|R| = n$.

%

The following lemma
and its generalisation (stated in Corollary~\ref{c:BipPW} in the next section) are
used frequently in later proofs.

\begin{lem}\lab{l:completeBipartite}
If a $2$-colouring of $K_C$ contains a monochromatic well-split $K_{2n^2,2n^2}$, then it contains a
monochromatic non-crossing copy of $L_{2n}$.
\end{lem}

\proof Let $A_1$ and $A_2$ be the two vertex parts of the monochromatic well-split
$K_{2n^2,2n^2}$. Without loss of generality, we assume that all the edges
between $A_1$ and $A_2$ are coloured blue.

We use an idea that was used to prove a quadratic upper bound on $R_g(C_n)$ and
other results on geometric Ramsey numbers~\cite{KPT,KPTV}.
We define partial orders $<_1$ on $A_1$ and $<_2$ on $A_2$ as follows.
A path $(p_i)_{i=1}^{\ell}$ on the vertices of $A_i$ is an \emph{increasing path}
if $p_1 \prec p_2 \prec \dots \prec p_{\ell}$.
Let $u <_i v$ for $u,v \in A_i$
if and only if there exists an increasing blue path starting in $u$ and ending at $v$.
Since $|A_i|=2n^2$ for $i=1,2$, by a lemma of Dilworth~\cite{Dilworth},
each of $(A_i,<_i)$ has either a chain on $n$ elements or an antichain on $2n$ elements.

By the definition of $(A_i,<_i)$,
any two vertices that are incomparable in $A_i$ are connected by a red edge.
Therefore, if $(A_i,<_i)$ contains an antichain with $2n$ elements,
then there exists a red convex complete geometric graph on $2n$ vertices.
By Lemma~\ref{l:gritzmannEtAl}, $K_C$ contains a red non-crossing copy of $L_{2n}$.

Thus, we may assume that none of  $(A_i,<_i)$ contains an antichain with $2n$ elements.
Then both  $(A_i,<_i)$ contain a chain with $n$ elements,
implying that each of $A_i$ contains an increasing blue path with $n$ vertices.
Let $(u_i)_{i=1}^{n}$ and $(v_i)_{i=1}^{n}$ be the increasing blue paths on $A_1$ and $A_2$, respectively.
These two paths together with the blue edges $\{u_{n+1-i},v_i\}$, $i\in [n]$
form a blue non-crossing copy of $L_{2n}$.
\qed \ss


\begin{lem}
\lab{l:longPath}
Let $N$ and $n$ be positive integers.
Let $G$ be the complete graph on a set $A$ of at least $nN$ vertices
and let $(A_i)_{i=1}^{n}$ be a partition of $A$ with $|A_i| \ge N$ for every $i \in [n]$.
Then for any 2-colouring of the edges of $G$,
either there is a red path $(u_i)_{i=1}^{n}$ with $u_i\in A_i$ for each $i \in [n]$
or for some $i \in [n-1]$ there exists a blue $K_{B_i, B_{i+1}}$ with
$B_i\subseteq A_i$, $B_{i+1}\subseteq A_{i+1}$ and $\min\{|B_i|, |B_{i+1}|\}\ge N/2$.
\end{lem}


\proof
Assume that there is no red path $(v_i)_{i=1}^{n}$ with $v_i \in A_i$ for each $i \in [n]$.
We call a vertex $v \in A_j$ \emph{good} if there is a red path $(v_i)_{i=1}^{j}$
with $v_j=v$ and $v_i \in A_i$ for every $i \in [j-1]$.
Every vertex in $A_1$ is good and all vertices in $A_n$ are bad.
Let $i$ be the largest integer such that at least half of the vertices of $A_i$ are good.
Then by the choice of $i$, at least half of the vertices of $A_{i+1}$ are bad.
Let $B_i$ denote the set of good vertices in $A_i$ and $B_{i+1}$ the set of bad vertices in $A_{i+1}$.
It follows that both $B_{i}$ and $B_{i+1}$ have size at least $N/2$
and all the edges between $B_{i}$ and $B_{i+1}$ are blue.
\qed \ss

\no{\bf Proof of Theorem~\ref{t:Ladder}}
Let $C$ denote a set of $32n^3$ points in convex position.
Arbitrarily choose a line that partitions $C$ into $C_1$ and $C_2$
each containing exactly $16n^3$ points.
Further, partition $C_1$ into $(A_i)_{i=1}^{2n}$ with $A_1\prec A_2\prec\cdots\prec A_{2n}$
and $|A_i|=8n^2$ for each $i \in [2n]$.
Partition $C_2$ into $(B_i)_{i=1}^{2n}$ with $B_{2n}\prec B_{2n-1}\prec\cdots\prec B_{1}$
and $|B_i|=8n^2$ for each $i \in [2n]$.

Colour each vertex $v\in A_i$ red if it is adjacent to at least half of the vertices in $B_i$ by a red edge. Otherwise, colour it blue.
We say that $A_i$ is red if at least half of the vertices in $A_i$ are coloured red.
Otherwise, we say that $A_i$ is blue.

Without loss of generality, at least half of the sets $A_i$ are red.
Let $(j_i)_{i=1}^{n}$ be an increasing sequence of indices such that each $A_{j_i}$ is red.
Let $D_i\subset A_{j_i}$ be the set of red vertices of $A_{j_i}$.
Thus $|D_i| \ge |A_{j_i}|/2 = 4n^2$ and for every vertex $v$ from $D_i$,
$|B_{j_i} \cap N_{red}(v)| \geq 4n^2$.

If for some $i \in [n-1]$, there exists a blue $K_{T_i, T_{i+1}}$
with $T_i\subseteq D_i$, $T_{i+1}\subseteq D_{i+1}$ and $\min\{|T_i|, |T_{i+1}|\}\ge 2n^2$,
then $K_C$ contains a blue non-crossing copy of $L_{2n}$ by Lemma~\ref{l:completeBipartite}.
Thus, by Lemma~\ref{l:longPath}, we can assume that we have a red path $(v_i)_{i=1}^{n}$ with
$v_i \in D_i$ for every $i \in [n]$.
So $v_1 \prec v_2 \prec \cdots \prec v_n$.
For each $i \in [n]$, let $F_i = B_{j_i} \cap N_{red}(v_i)$.
So $|F_i| \geq  4n^2$.
If there exists a blue $K_{T_i, T_{i+1}}$ with $T_i\subseteq F_i$, $T_{i+1}\subseteq F_{i+1}$
and $\min\{|T_i|, |T_{i+1}|\}\ge 2n^2$, then the proof is complete by Lemma~\ref{l:completeBipartite}.
Thus by Lemma~\ref{l:longPath}, we only need to consider the case when
there is a red path $(w_i)_{i=1}^{n}$ with $w_i \in F_i$ for every $i \in [n]$.
We have $w_n \prec w_{n-1} \prec \cdots w_1$ and so the two paths $(v_i)_{i=1}^{n}$ and $(w_i)_{i=1}^{n}$
together with the edges $\{v_i,w_i\}$, $i \in [n]$ form a red non-crossing copy of $L_{2n}$.
\qed

\subsection{Ordered Ramsey theory} 
\lab{s:ladder:ordered}

The proof of Theorem~\ref{t:Ladder} shows that the ordered Ramsey number $R_o(L_{2n}, L_{2n})$ of the ladder graph $L_{2n}$ with vertices ordered 
$v_1 \prec v_2 \prec \dots \prec v_n \prec u_n \prec u_{n-1} \prec \dots \prec u_1$ is at most $32 n^3$.
The ideas used in the proof of Theorem~\ref{t:Ladder} can be applied to give a subexponential upper bound 
on the ordered Ramsey numbers of arbitrarily ordered paths. 
\begin{thm}
\label{t:OrderedRamsey}
Let $K_n$ be the ordered complete graph on $n$ vertices and let $P_m$ be an arbitrarily ordered path on $m$ vertices.
Then $R_o(K_n, P_m) \le 2^{\lceil \log_2(n) \rceil \cdot (\lceil \log_2(m) \rceil + 1)}$.
\end{thm}

\proof
If $n \le 2$, the claim holds trivially. We proceed by induction on $n$ while $m$ remains fixed.

It is enough to show that  $R_o(K_n, P_m) \le 2^{\log_2(n) \cdot (\log_2(m) + 1)}$
for values $n$ and $m$ of the form $n = 2^k$ and $m = 2^{\ell}$ for some integers $k$ and $\ell$.
Let $R = 2^{k\cdot(\ell+1)}$.
Let $K_R$ be the complete ordered graph on $R$ vertices with $2$-coloured edges.
We split the vertices of $K_R$ into $m$ intervals $V_1$, \ldots, $V_m$, 
each containing $2^{(k-1)\cdot(\ell+1)+1}$ consecutive vertices.
Let $p_1 \prec p_2 \prec \ldots \prec p_m$ be the vertices of $P_m$.
Then the edges of $P_m$ are $\{p_{\pi(1)}, p_{\pi(2)}\}, \{p_{\pi(2)}, p_{\pi(3)}\}, \ldots, \{p_{\pi(m-1)}, p_{\pi(m)}\}$ 
for some permutation $\pi:[m]\rightarrow [m]$.
We let $A_i = V_{\pi(i)}$ for every $i \in [n]$.
By Lemma~\ref{l:longPath}, we either find a red copy of $P_m$, in which case the claim holds,
or we find a pair of intervals $A_{i}$, $A_{i+1}$ satisfying the following.
There are sets $L \subset A_{i}$ and $R \subset A_{i+1}$
of size $|L|, |R| \ge 2^{(k-1)\cdot(\ell+1)}$ such that all the edges between $L$ and $R$ are blue.

By the induction hypothesis, $R_o(K_{n/2}, P_{m}) \le 2^{(k-1)\cdot(\ell+1)}$.
Thus in each of $L$ and $R$, we either find a red copy of $P_m$, or a blue copy of $K_{n/2}$.
If either $L$ or $R$ contains a red copy of $P_m$, the claim holds.
Otherwise both $L$ and $R$ contain a blue copy of $K_{n/2}$ and so $L \cup R$ is a blue copy of $K_n$.
\qed

\begin{cor}
\label{c:OrderedRamsey}
Let $P_n$ be a path on $n$ arbitrarily ordered vertices.
Then $R_o(P_n, P_n) \le 2^{\lceil \log_2(n) \rceil \cdot (\lceil \log_2(n) \rceil + 1)}$.
\end{cor}

\subsection{General geometric position}
\lab{s:ladder:general}


\begin{thm}\lab{gent:Ladder}
The geometric Ramsey number of the ladder graph $L_{2n}$ satisfies $R_g(L_{2n})=O(n^{10})$.
\end{thm}

\begin{defn}
Two sets of points $A$ and $B$ in the plane are {\em mutually avoiding} if $|A|,|B| \geq 2$
and no line subtended by a pair of points in $A$ intersects the convex hull of $B$, and vice versa.
See Fig.~\ref{fig:muavsets}.
\end{defn}

\onefigure{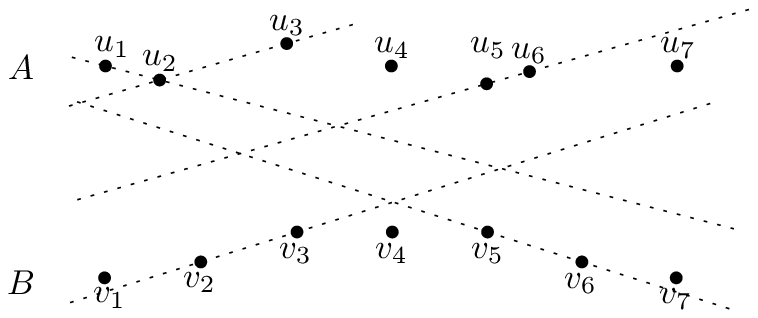}{1}{An example of mutually avoiding sets $A$ and $B$. 
Some lines subtended by pairs of points from $A$ and pairs of points from $B$ are shown. }{fig:muavsets}

A simple example of a pair of mutually avoiding sets are sets $A$ and $B$ such that $A\cup B$ is
in convex position and $A$ and $B$ can be separated by a straight line.

Observe that for any mutually avoiding pair $(A,B)$, every point in $A$ ``sees'' all the vertices in $B$ in
the same order and vice versa.
That is, there are unique total orders $u_1\prec u_2 \prec \dots \prec u_{|A|}$ of the points in $A$
and $v_1\prec v_2 \prec \dots \prec v_{|B|}$ of the the points in $B$ such that every point in $B$ ``sees'' $u_1,\ldots,u_{|A|}$ consecutively in a clockwise order before seeing any vertex in $B$, whereas every point in
$A$ ``sees'' $v_1,\ldots,v_{|B|}$ consecutively in a counterclockwise order before seeing any vertex in $A$.
A path $(p_i)_{i=1}^{\ell}$ in either $A$ or $B$ is an \emph{increasing path}
if $p_1 \prec p_2 \prec \cdots \prec p_{\ell}$.

For any two sets of vertices $A_1,A_2$ both contained in $A$ (or $B$), we write
$A_1\prec A_2$ if and only if for every $u\in A_1$ and $v\in A_2$, $u\prec v$.
Let $U$ be $A$ or $B$.
A sequence $(U_i)_{i=1}^{k}$ of subsets of $U$ is an {\em increasing sequence}
if $U_1\prec \cdots \prec U_k$.
An increasing sequence $(U_i)_{i=1}^{k}$ of subsets of $U$ is an {\em increasing partition} of $U$
if $\bigcup_{i \in [k]} U_i=U$.

The following proposition follows from the definition of a pair of mutually avoiding sets.
\begin{prop} \lab{prop:muav} {\ }
Assume $A$ and $B$ are mutually avoiding. Then
\begin{enumerate}
\item
An increasing path $(p_i)_{i=1}^{\ell}$ does not cross itself.
\item
Let $u, u'\in A$ with $u \prec u'$ and let $v, v'\in B$.
Then the two edges $\{u, v\}$ and $\{u', v'\}$ cross if and only if $v' \prec v$.
\item
Let $u, u'\in A$ with $u \prec u'$ and let $w\in A$ such that
$w \prec u$ or $u' \prec w$ or $w\in \{u,u'\}$.
Let $v\in B$.
Then the two edges $\{u, u'\}$ and $\{w, v\}$ do not cross.
\end{enumerate}
\end{prop}

The following corollary follows directly from Proposition~\ref{prop:muav}.
\begin{cor} \lab{cor:muavLadder}
Assume $A$ and $B$ are mutually avoiding. Let $P_u = (x_{i})_{i=1}^{n}$ be an increasing path in $A$
and let $P_v = (y_{i})_{i=1}^{n}$ be an increasing path in $B$.
Then the ladder graph composed of the paths $P_u$ and $P_v$ and edges $\{\{x_{i},y_{i}\}: i\in [n]\}$ is non-crossing.
\end{cor}

Given a set of points in general position, the following theorem guarantees the existence of two mutually avoiding subsets of relatively large sizes.

\begin{thm}[Aronov et al.~1994~\cite{AronovEtAl}]\lab{t:avoiding}
Let $A'$ and $B'$ be two sets of points separated by a line,
each of size $6 n^2$. Then there exist mutually avoiding sets
$A\subset A'$ and $B\subset B'$ such that $A$ and $B$ are both of
size $n$.
\end{thm}

An embedding of the complete bipartite graph $K_{m,n}$ on a set of points in general position 
is {\em well-split} if the two sets of points representing the two vertex parts 
are mutually avoiding.
By Lemma~\ref{l:gritzmannEtAl} and Corollary~\ref{cor:muavLadder}, we have the following
generalisation of Lemma~\ref{l:completeBipartite}.
\begin{lem}\lab{genl:completeBipartite}
If a $2$-colouring of $K_P$ contains a monochromatic well-split $K_{2n^2,2n^2}$, then it contains a
monochromatic non-crossing $L_{2n}$.
\end{lem}
\proof
Let $A_1$ and $A_2$ be the two vertex parts of the monochromatic well-split $K_{2n^2,2n^2}$.
Without loss of generality, assume $K_{A_1,A_2}$ is blue.
By applying the Dilworth's lemma~\cite{Dilworth} in the same way as in the proof of
Lemma~\ref{l:completeBipartite}, we either find a red $K_{2n}$
or blue increasing paths $(u_i)_{i=1}^{n}$ in $A_1$ and $(v_i)_{i=1}^{n}$ in $A_2$.
In the first case, we get a red $L_{2n}$ by Lemma~\ref{l:gritzmannEtAl}
and in the second case a blue $L_{2n}$ by Corollary~\ref{cor:muavLadder}.
\qed \ss

A complete geometric bipartite graph $K_{L,R}$ is said to be \emph{separable}
if $L$ and $R$ can be separated by a line.
Notice that if $L \cup R$ is in convex position, then $K_{L,R}$ is separable if and only if it is well-split.
Obviously, every complete bipartite geometric graph $K_{L,R}$ contains a separable complete bipartite graph
with parts of sizes $|L|/2$ and $|R|/2$.
However, all complete bipartite geometric graphs that we encounter in subsequent proofs are separable,
so we state the following corollary of Theorem~\ref{t:avoiding} and Lemma~\ref{genl:completeBipartite} for
separable complete bipartite graphs only.

\begin{cor}\lab{c:completeBipartite}
Every $2$-colouring of $K_P$ containing a monochromatic separable $K_{24n^4,24n^4}$
contains a monochromatic non-crossing copy of $L_{2n}$.
\end{cor}

\medskip

\no{\bf Proof of Theorem~\ref{gent:Ladder}.\ }
Let $G$ be the complete geometric graph on vertex set $P$,
where $P$ is a set of $cn^{10}$ points in general position,
where $c$ is some sufficiently large absolute constant.
By Theorem~\ref{t:avoiding}, there exist two subsets $S_u,S_v\subset
P$, such that $S_u$ and $S_v$ are mutually avoiding and
$|S_u|=|S_v|=c_1n^5$ for some $c_1\ge \sqrt{c/6}$.

The proof of Theorem~\ref{gent:Ladder} is analogous to that of Theorem~\ref{t:Ladder}
with $S_u$ and $S_v$ having the role of $C_1$ and $C_2$.
Let $(A_i)_{i=1}^{2n}$ be the increasing partition of $S_u$ with $|A_i|=c_1n^4/2$ for each $i \in [2n]$.
Let $(B_i)_{i=1}^{2n}$ be the increasing partition of $S_v$ with $|B_i|=c_1n^4/2$ for each $i \in [2n]$.
As in the proof of Theorem~\ref{t:Ladder}, we find one colour, that we assume to be red,
an increasing sequence $(D_i)_{i=1}^{n}$ of subsets of $S_v$
and an increasing sequence $(j_i)_{i=1}^{n}$ of integers from $[2n]$ satisfying the following.
For every $i \in [n]$, $|D_i| \ge c_1 n^4/4$
and every vertex in $D_i$ is adjacent to at least half of the vertices of $B_{j_i}$ by a red edge.
By Lemma~\ref{l:longPath}
either there is a blue copy of $K_{c_2 n^4,c_2 n^4}$, where $c_2\ge c_1/8\ge \sqrt{c}/(16\sqrt{3})$
or there is a red path $(v_i)_{i=1}^{n}$ with $v_i \in D_i$ for every $i \in [n]$.
In the first case the proof is complete by Corollary~\ref{c:completeBipartite}.
In the second case we let $T_i = B_{j_i} \cap N_{red}(v_i)$.
Now we apply Lemma~\ref{l:longPath} on $(T_i)_{i=1}^n$ and either find a blue $K_{c_2 n^4,c_2 n^4}$
or a red path $(w_i)_{i=1}^n$ with $w_i \in T_i$ for every $i \in [n]$.
In the first case the proof is complete by Corollary~\ref{c:completeBipartite}
and in the second by Corollary~\ref{cor:muavLadder} using $(v_i)_{i=1}^{n}$ and $(w_i)_{i=1}^{n}$
as the two paths of the ladder graph.
 \qed

\section{Generalisation to pathwidth-2 outerplanar triangulations}\lab{s:generalisation}

An \emph{outerplanar triangulation} $G$ is a planar graph that can be drawn in the plane in such a way that the outer face is incident with all the vertices of $G$ and every other face is incident with exactly three vertices.

The pathwidth of a graph was first defined by Robertson and Seymour~\cite{RS} as follows.
A \emph{path decomposition} of a graph $G$ is a sequence $(G_i)_{i=1}^{m}$ of subgraphs of $G$
such that each edge of $G$ is in at least one of $G_i$
and for every vertex $v$ of $G$, the set of graphs $G_i$ containing $v$ forms a contiguous subsequence of
$(G_i)_{i=1}^{m}$.
The \emph{pathwidth} of a graph $G$ is the smallest $k$ such that $G$ has a path decomposition in which
every $G_i$ has at most $k+1$ vertices. Let $pw(G)$ denote the pathwidth of $G$.
A \emph{pathwidth-$k$} graph is a graph of pathwidth at most $k$.

For every $k$, the class of graphs of pathwidth at most $k$ is a minor-closed class.
Every such class can be characterised by a finite list of forbidden minors 
by the graph minor theorem~\cite{RSXX}.
A characterisation of the class of pathwidth-$2$ graphs with $110$ forbidden minors 
was provided by Kinnersley and Langston~\cite{KL94}.


Simplified characterisations of pathwidth-$2$ graphs were obtained 
recently by Bar\'{a}t, Hajnal, Lin and Young~\cite{BHLY} and Bir{\'o}, Keller and Young~\cite{BKY}.
We use these characterisations to provide an equivalent definition of 
pathwidth-$2$ outerplanar triangulations that will be used in our proofs.

\begin{defn} \label{def:PW2}
Let $\PW(n)$ be the class of outerplanar triangulations $G$ on $n$ vertices whose vertices 
can be decomposed into two disjoint sets $V_u \cup V_v = V(G)$ such that 
the subgraphs induced by the two sets, $P_u = G[V_u]$ and $P_v = G[V_v]$, are paths. 
See an example in Fig.~\ref{fig:expw}.
\end{defn}

\onefigure{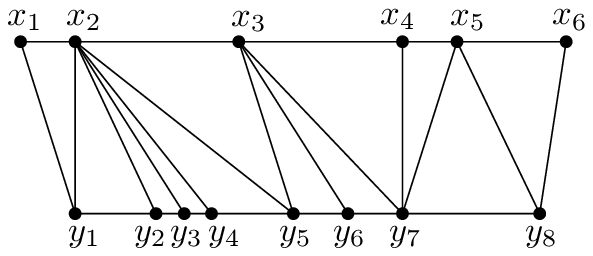}{1}{An example of a pathwidth-$2$ outerplanar triangulation.}{fig:expw}

\begin{prop}\lab{p:pathwidth}
A graph $G$ is a pathwidth-$2$ outerplanar triangulation if and only if $G \in \PW(n)$.
\end{prop}
\proof
A \emph{track} is a graph composed of two \emph{rails} and several \emph{cross-ties}.
The two rails are paths $(x_i)_{i=1}^{n_1}$ and $(y_i)_{i=1}^{n_2}$.
A cross-tie is a path of length one or two that connects $x_i$ with $y_j$, for some 
$i \in [n_1]$ and $j \in [n_2]$.
The cross-ties further satisfy that for every $i, i', j, j'$ with $i < i'$, whenever 
one cross-tie connects $x_i$ to $y_j$ and another connects $x_{i'}$ to $y_{j'}$,
then $j\le j'$.
The middle vertex of a cross-tie of length two has no neighbours other than $x_i$ and $y_j$.
Additionally, there always is a cross-tie of length one connecting $x_1$ to $y_1$ 
and another connecting $x_{n_1}$ to $y_{n_2}$.
Bar\'{a}t et~al.~\cite{BHLY} prove that a graph is a $2$-connected pathwidth-$2$ graph 
if and only if it is a track.

Notice that every outerplanar triangulation is Hamiltonian and thus $2$-connected.
Observe also that every $G \in \PW(n)$ satisfies the definition of a track.
It remains to show that if a track $G$ on $n$ vertices is an outerplanar triangulation
then $G \in \PW(n)$.

If the track $G$ has a cross-tie of length two between $x_1$ and $y_1$, then the middle 
vertex $t$ of the cross-tie can be added to one of the tracks to form, for example, 
the track $t, x_1, \ldots, x_{n_1}$.
Thus we can assume, that the track $G$ has no cross-tie of length two connecting 
$x_1$ to $y_1$ or $x_{n_1}$ to $y_{n_2}$.
Then, since $G$ is outerplanar, it has no cross-tie of length two.
Therefore the outerplanar triangulation $G$ satisfies the definition of graphs 
from $\PW(n)$.
\qed\ss

The following is a corollary of Property~\ref{prop:muav}.

\begin{cor}\lab{cor:muavPW2}
Let $G\in \PW(n)$ and let $G$ be composed of induced paths $P_u=(x_i)_{i=1}^{n_1}$,
$P_v=(y_i)_{i=1}^{n_2}$ and edges between vertices of $P_u$ and vertices of $P_v$.
Let $(A, B)$ be a pair of mutually avoiding sets.
Let $(u_{i})_{i=1}^{n_1}$ be an increasing path in $A$
and $(v_{i})_{i=1}^{n_2}$ an increasing path in $B$.
Then by mapping every $x_i$ on $u_i$ and every $y_i$ on $v_i$ we obtain a non-crossing embedding of $G$.
\end{cor}

By Corollary~\ref{cor:muavPW2},
Lemma~\ref{l:completeBipartite} generalises to an arbitrary graph $G \in \PW(n)$.
\begin{cor}\lab{c:wellSplitBipPW}
If $K_P$ with $2$-coloured edges contains a monochromatic well-split $K_{n^2,n^2}$,
then it contains a monochromatic non-crossing copy of every $G$ from $\PW(n)$.
\end{cor}

Then, by Theorem~\ref{t:avoiding}, we also generalise Corollary~\ref{c:completeBipartite}.
\begin{cor}\lab{c:BipPW}
If $K_P$ with $2$-coloured edges contains a monochromatic separable $K_{6n^4,6n^4}$,
then it contains a monochromatic non-crossing copy of every $G$ from $\PW(n)$.
\end{cor}

We obtain upper bounds for the geometric Ramsey numbers of graphs $G\in\PW(n)$, both in the convex case and in the general case. These two upper bounds follow directly from the following key lemma.

\begin{lem}\lab{l:pathwidthboth}
Let $G$ be a subgraph of a graph $G' \in \PW(n)$.
Let $m \geq n^2$ and let $S_u$ and $S_v$ be two mutually avoiding
sets of $10 m^2 n^3$ points each. Then every $2$-colouring of
the complete geometric graph on $S_u \cup S_v$ either contains a monochromatic
$G$ or a monochromatic separable $K_{m,m}$.
\end{lem}

We leave the technical proof of Lemma~\ref{l:pathwidthboth} to the next section.

\begin{thm}\lab{t:pathwidthconvex}
For any $G\subseteq G'\in \PW(n)$, $R_c(G)\le 20 n^{7}$.
\end{thm}

\proof
Let $S$ be a set of $20 n^7$ points in convex position. We cut the set
$S$ by a line into sets $S_u$ and $S_v$ of size $10 n^7$ each.
Then $S_u$ and $S_v$ are mutually avoiding.
Moreover if either $S_u$ or $S_v$ contains a monochromatic separable and thus well-split $K_{n^2,n^2}$,
then $S$ contains a monochromatic non-crossing $G$ by Corollary~\ref{c:wellSplitBipPW}.

Therefore, by Lemma~\ref{l:pathwidthboth} with $m=n^2$, $S$ contains a monochromatic
non-crossing copy of $G$.
\qed\ss

\begin{thm}\lab{t:pathwidthgeneral}
For any $G\subseteq G'\in \PW(n)$, $R_g(G)\le O(n^{22})$.
\end{thm}

\proof
Let $S$ be a set of $10^2 6^5 n^{22}$ points in general position. By Theorem~\ref{t:avoiding},
$S$ contains mutually avoiding sets $S_u$ and $S_v$ of size $10\cdot 6^2n^{11}$ each.
If $S$ contains a monochromatic separable $K_{6n^4, 6n^4}$, then it contains
a monochromatic non-crossing copy of $G$ by Corollary~\ref{c:BipPW}.

Therefore, by Lemma~\ref{l:pathwidthboth} with $m = 6 n^4$, $S$ contains a monochromatic
non-crossing copy of $G$.
\qed\ss

\begin{rem*}
Notice that not every pathwidth-$2$ outerplanar graph is a subgraph of a 
pathwidth-$2$ outerplanar triangulation. 
See Fig.~\ref{fig:notTriang}.
\end{rem*}

\onefigure{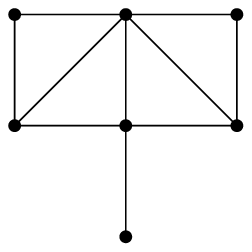}{1}{A pathwidth-$2$ outerplanar graph that is not a subgraph of a 
pathwidth-$2$ outerplanar triangulation. }{fig:notTriang}

\section{Proof of Lemma~\ref{l:pathwidthboth}}

Without loss of generality, we only need to consider $G\in\PW(n)$.

If $n \le 3$ then the result is trivial, so we assume $n \ge 4$, which will be used in
several estimates.

Let $(u_i)_{i=1}^{\ell-1}$ and $(v_i)_{i=1}^{\ell}$ be sequences of vertices such that
$u_i \in P_u$ and $v_i \in P_v$ for every $i \in [\ell]$, $v_1 = y_1$, $v_{\ell} = y_{n_2}$
and the alternating sequence of vertices $v_1, u_1, v_2, u_2, \dots, u_{\ell-1}, v_{\ell}$
forms a path in $G$.
Note that $2\ell-1 \leq n$ and that such a path is unique
and can be constructed by starting at $y_1$
and always continuing to the largest neighbour in the other of the sets $P_u, P_v$
until reaching $y_{n_2}$.
Vertices $u_1,\dots u_{\ell-1}, v_1,\dots v_{\ell}$
are called the \emph{stem vertices} and all the other vertices of $G$ are the
\emph{leaf vertices}.
Notice that every leaf vertex of $P_u$ has exactly one
neighbour in $P_v$ and vice versa.

\onefigure{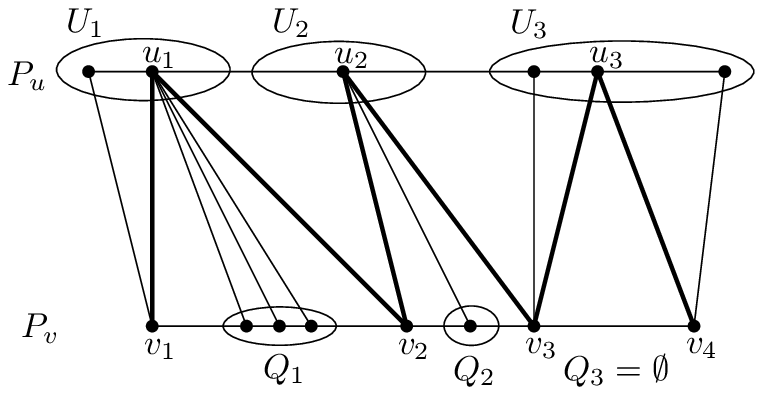}{1}{Labelling and grouping of the vertices of the graph from Fig.~\ref{fig:expw}, which is in $\PW(14)$.
The path connecting stem vertices is represented by a heavier line. }{fig:pw-start}

Refer to Fig.~\ref{fig:pw-start}.
We cut $P_u$ into a sequence of subpaths $(U_i)_{i=1}^{\ell-1}$, where $U_1$ contains
$u_1$ and all the vertices preceding $u_1$ on $P_u$, $U_{\ell-1}$ contains all the vertices
after $u_{\ell-2}$ and for $i \in [2,\ell-2]$, $U_i$ contains $u_{i}$ and the vertices
strictly between $u_{i-1}$ and $u_i$.
Let $(Q_i)_{i=1}^{\ell-1}$ be the sequence of subpaths created by removing
vertices $v_i$ from $P_v$. That is, for every $i \in [\ell-1]$, $Q_i$ contains
the leaf vertices in between $v_i$ and $v_{i+1}$.
Let $f_i = |Q_i|$.

\begin{cl}\lab{cl:initialPartition}
There exists a colour $c\in\{\mbox{blue,red}\}$
and sequences  of vertex sets $A'_1, \dots, A'_{\ell} \subset S_v$,
$M_1, \dots, M_{\ell-1} \subset S_v$ and
$B_1, \dots, B_{\ell-1} \subset S_u$ with
\begin{align*}
 & B_1\prec B_2\prec\cdots\prec B_{\ell-1} \\
 & A'_{1}\prec M_{1}\prec A'_{2}\prec\cdots\prec A'_{\ell-1}\prec M_{\ell-1} \prec A'_{\ell}
\end{align*}
that satisfy the following conditions.
\begin{enumerate}
\item $\forall i \in [\ell-1]: |B_i| = 8 mn^3$;
\item $\forall i \in [\ell]: |A'_i| = 4 mn^2$;
\item $\forall i \in [\ell-1]: |M_i| = 9 m^2n^2$;
\item \label{it:halfneigh} $\forall i \in [\ell-1] :\forall v \in A'_i:$ $|N_{c}(v) \cap B_{i}| \ge |B_{i}|/2$.
\end{enumerate}
\end{cl}

\proof
Let $(Z_k)_{k=1}^{2\ell-1}$ be the increasing partition of $S_u$ with parts of size $8 mn^3$.
Then we take sequences $(D_k)_{k=1}^{2\ell-1}$, $(C_k)_{k=1}^{2\ell-2}$ of subsets of $S_v$ satisfying
$D_{1}\prec C_{1}\prec D_{2}\prec\cdots\prec D_{2\ell-2}\prec C_{2\ell-2} \prec D_{2\ell-1}$ with
$|D_k| = 8 mn^2$ for every $k \in [2\ell-1]$ and $|C_k| = 9 m^2n^2$ whenever $k \in [2\ell-2]$.

The \emph{colour of a vertex $v \in D_k$} is the colour of the majority
of the edges between $v$ and the vertices of $Z_k$ and it is red in case of a tie.
The \emph{colour of $D_k$} is the colour of the majority of the vertices $v \in D_k$
and it is red in case of a tie.

We fix $c$ to be the colour such that at least half of the sets $D_k$ have colour $c$.
Let $(k_i)_{i=1}^{\ell}$ be an increasing sequence of indices such that
for every $i \in [\ell]$, $D_{k_i}$ has colour $c$.
For each $i \in [\ell]$ let $A'_i$ be the set of vertices of $D_{k_i}$ with colour $c$.
Then $|A'_i| \ge 4 mn^2$.
Let $B_{i} = Z_{k_i}$ for every $i \in [\ell-1]$ and $M_{i} = C_{k_i}$ for every $i \in [\ell-1]$.
It is easy to verify that the sets $A'_i$, $B_i$ and $M_i$ satisfy the requirements.
\qed \ss

\begin{obs}\lab{o:intersectingSets} Let $N$, $k$ and $t$ be positive
integers. Let $S$ be a set of size $N$ and let $T_1,\ldots,T_k$ be
sets such that $|S\cap T_i|\le t$ for every $i\in [k]$. Then
\[
|S \setminus \bigcup_{i=1}^{k} T_i| = |\bigcap_{i=1}^{k} (S \setminus T_i)| \ge N-tk.
\]
\end{obs}

The observation is applied several times for some colour $c$ and a set $V$ of vertices
in the following way.
We set $T_i = N_c(v_i)$ where $\{v_1, \dots, v_k\}$ are the vertices of $V$ with
the fewest $c$-neighbours in $S$.
The observation says that if $k$ is large and every $v_i$ has few $c$-neighbours in $S$,
then we find a large complete bipartite graph in the other colour.

Without loss of generality, we assume that Claim~\ref{cl:initialPartition} holds with $c=\mbox{red}$.
Let $(A'_i)_{i=1}^{\ell}$, $(B_i)_{i=1}^{\ell}$ and $(M_i)_{i=1}^{\ell-1}$
be the sequences that satisfy the conditions of the claim.

The rest of the proof proceeds in several phases.
In each phase we either immediately find a blue well-split $K_{n^2, n^2}$ implying a monochromatic $G$,
a blue separable $K_{m,m}$,
or we move closer to finding a non-crossing embedding $\phi: V(G) \rightarrow S_u \cup S_v$ of $G$ with all edges red.
The mapping $\phi$ maps $v_i$ on some point of $A'_i$ for each $i \in [\ell]$ and for each $i \in [\ell-1]$,
the vertices of $U_i$ are mapped on some points of $B_i$ and vertices of $Q_i$ on some points of $M_i$.
In some phases, the embedding of some vertices of $G$ is selected
and this will then remain fixed for the rest of the proof.

\begin{cl}\lab{cl:finerPartition}
Either the complete geometric graph on $S_u \cup S_v$ contains a monochromatic noncrossing $G$
or there is a sequence of sets $(A_i)_{i=1}^{\ell}$ with $A_i\subseteq A'_{i}$
for every $i \in [\ell]$ that satisfies the following conditions.
\begin{enumerate}
\item $\forall i: |A_i| = 2m$;

\item \label{it:inters}$\forall i \in [\ell-1], \forall u \in A_i, \forall v \in A_{i+1}:$
$u$ and $v$ have at least $3nm$ common red neighbours in $B_{i}$;

\end{enumerate}
\end{cl}

\proof
To find the sets $A_i$, we proceed in $\ell$ steps, unless we find a red $G$ earlier.

In the first step, we let the set $A_{1}$ be an arbitrary subset of $A'_{1}$ of size $2m$.

At the beginning of step $j$, $j > 1$, we have sets $A_{1}, \dots, A_{j-1}$ each of size $2m$
and such that the requirement~\ref{it:inters} of the claim is satisfied for all $i < j-1$.
A vertex $w \in A'_{j}$ is \emph{compatible} with $v \in A_{j-1}$ if $u$ and $v$ have
at least $3nm$ common red neighbours in $B_{j}$.
We distinguish two cases.

The first case occurs when there is a vertex $v \in A_{j-1}$
and a set $W = \{w_1, \dots, w_{n^2}\}$ of vertices of $A'_j$ incompatible with $v$.
Let $S = N_{red}(v) \cap B_{j}$ and for every $i\in[n^2]$, let $T_i = N_{red}(w_i) \cap B_{j}$.
Let $C = S \setminus \bigcup_{i=1}^{n^2} T_i$.
Since the vertices of $W$ are incompatible with $v$, we can apply Observation~\ref{o:intersectingSets}
on $S$ and $\{T_1, \dots, T_{n^2}\}$ with $t = 3nm$, $N = 4mn^3$ and $k = n^2$
and obtain $|C| \ge 4mn^3-3mn^3 \ge n^2$.
All edges between $W$ and $C$ are blue, thus $K_{W,C}$ forms a blue well-split $K_{n^2,n^2}$
and so $K_{S_u \cup S_v}$ contains a monochromatic noncrossing $G$ by Corollary~\ref{c:wellSplitBipPW}.

In the second case, for every vertex $u \in A_{j-1}$ at most $n^2$ vertices of $A'_j$ are incompatible.
Thus the number of vertices of $A'_j$ compatible with every $u \in A_{j-1}$
is at least $4mn^2 - 2m n^2 \geq 2m$.
We can thus let $A_j$ be the set of some $2m$ vertices of $A'_j$ compatible with every $v \in A_{j-1}$.
\qed \ss

Let $(A_i)_{i=1}^{\ell}$ be the sequence of sets satisfying the conditions of
Claim~\ref{cl:finerPartition}.

To provide an exposition of the rest of the proof,
we first prove Lemma~\ref{l:pathwidthboth} for the case when there is no leaf vertex on $P_v$.

\begin{cl}\lab{cl:noLeaf} Assume $P_v$ contains no leaf vertices, then there exists a blue separable $K_{m,m}$
or a monochromatic non-crossing $G$.
\end{cl}

\proof
We assume that neither $S_u$ nor $S_v$ contains a blue separable $K_{m,m}$.
By Lemma~\ref{l:longPath}, we find a red path $(a_i)_{i=1}^{\ell}$, where each $a_i \in A_i$.
Then for every $i \in [\ell-1]$ we take the set $R_i \subseteq B_i$
of $3nm$ common red neighbours of $a_i$ and $a_{i+1}$.
For every $i \in [\ell-1]$, let $\mathcal{R}_i$ be an increasing partition of $R_i$
with $|U_i|$ parts of size at least $2m$ each.
By Lemma~\ref{l:longPath} we find an increasing red path $(r_i)_{i=1}^{|P_u|}$
with exactly one vertex in each set in $\bigcup_{i=1}^{\ell} \mathcal{R}_i$.
Then we map every $u_i$ on $r_i$ and every $x_i$ on $a_i$ to obtain a red copy of $G$.
\qed\ss

The rest of this section deals with the leaf vertices on $P_v$.

For each $i \in [\ell-1]$ such that $f_i>1$, we take an increasing partition $(M_{i,j})_{j=1}^{f_i}$ of $M_i$
with $|M_{i,1}|, |M_{i,f_i}| \ge 4m^2 n^2$ and $|M_{i,j}| \ge 3 mn^2$ for every $j \in [2,f_i-1]$.

Let $\gamma$ be the colouring of the edges of $K_{S_u \cup S_v}$.
We define a new edge colouring $\gamma'$ of the edges of $K_{S_u \cup S_v}$
according to the following cases.
\begin{enumerate}
\item The edge $e$ connects a vertex $v \in A_i$ and a vertex $w \in A_{i+1}$ such that $f_i = 1$.
We colour $e$ red if and only if $v$ and $w$ have at least $n^2$ common
red neighbours in $M_i$.
\item The edge $e$ connects a vertex $v \in A_i$ and a vertex $w \in A_{i+1}$ such that $f_i \ge 2$.
We colour $e$ red if and only if $|N_{red}(v) \cap M_{i,1}| \ge 3 mn^2$
and $|N_{red}(w) \cap M_{i,f_i}| \ge 3 mn^2$.
\item Otherwise $\gamma'(e) = \gamma(e)$.
\end{enumerate}

\begin{cl}
\lab{cl:atleast1leaf}
If there exist sets $L\subseteq A_i$ and $R\subseteq A_{i+1}$
with $|L|=|R|=m$ and all edges between $L$ and $R$ blue under $\gamma'$,
then there exists a blue separable $K_{m,m}$ in $S_v$ under $\gamma$.
\end{cl}
\proof
We distinguish three cases.
\begin{enumerate}
\item
We have $f_i = 1$.
If every $v \in L$ has fewer than $2 n^2m$ red neighbours in $M_i$ under $\gamma$,
then there are at least $9m^2n^2-2n^2m\cdot m\ge m$ vertices in $M_i$ that are
connected by blue edges to every vertex in $L$ under $\gamma$.
This implies the existence of a blue separable $K_{m,m}$ in $S_v$.
Otherwise, there exists a vertex $v \in L$ with at least $2 n^2 m$ red neighbours in $M_i$.
Let $N\subseteq M_i$ denote the set of these neighbours of $v$.
Since every edge between $L$ and $R$ is blue under $\gamma'$,
each $w \in R$ is connected by red edges to at most $n^2$ vertices in $N$.
Thus there are at least $2n^2m-n^2\cdot m\ge m$ vertices of $N$
that are connected by blue edges to each vertex in $R$.
Thus we have a blue separable $K_{m,m}$ in $S_v$ under $\gamma$.
\item
We have $f_i \ge 2$.
Either each point of $L$ has fewer than $3mn^2$ red neighbours in $M_{i,1}$
or each point of $R$ has fewer than $3 mn^2$ red neighbours in $M_{i,f(i)}$.
Without loss of generality, the first case occurs and
then there are at least $4m^2n^2-3mn^2\cdot m\ge m$ points in
$M_{i,1}$ connected by blue edges to every point of $L$.  
\item
We have $f_i = 0$.
Then $\gamma'$ is equal to $\gamma$ on all the edges between $A_i$ and $A_{i+1}$.
This implies the existence of a blue separable $K_{m,m}$. 
\qed
\end{enumerate}

By Claim~\ref{cl:atleast1leaf}, we may assume that
under $\gamma'$, there exists no $i \in [\ell]$ for which some two sets $L\subseteq
A_i$ and $R\subseteq A_{i+1}$ with $|L|=|R|=m$ would form a blue $K_{m,m}$.
Then by Lemma~\ref{l:longPath}, we can map each vertex $v_i$ on some point $\phi(v_i) \in A_i$ 
in such a way that $(\phi(v_i))_{i=1}^{\ell}$ is a red path under $\gamma'$.
For every $i \in [\ell-1]$, let $H'_i \subseteq B_i$ be a set of $3 nm$ common red
neighbours of $\phi(v_{i})$ and $\phi(v_{i+1})$
and let $H'_{\ell} = B_{\ell} \cap N_{red}(\phi(v_{\ell}))$.

In what follows, for each $i \in [\ell-1]$ we define a vertex set
$H_i\subseteq H'_i$ in which we then embed $U_i$.
If $f_i=1$ we define a vertex set $\widetilde{M}_i \subset M_i$ in which
we embed the only vertex of $Q_i$.
If $f_i\ge 2$ we define a sequence of vertex sets
$\widehat{M}_{i,1},\ldots,\widehat{M}_{i,f_i} \subset M_i$, and on each one of these sets,
we embed one of the leaf vertices from $Q_i$.
Refer to Fig.~\ref{fig:pw-sets}.

\onefigure{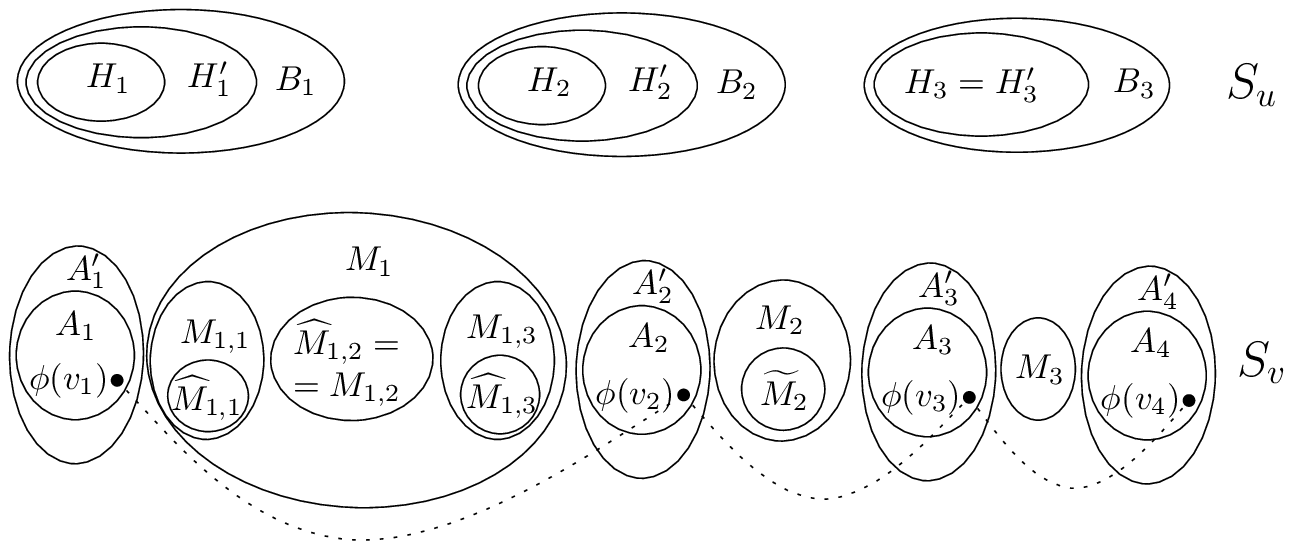}{1}{Embedded vertices $v_1, \ldots, v_4$ of the graph from Fig.~\ref{fig:pw-start}. 
This embedding fixes sets $\widetilde{M}_i$, $\widehat{M}_{i,j}$ and $H_i$. 
The dashed lines form a red path in the colouring $\gamma'$.}{fig:pw-sets}

If $f_i=0$, we let $H_i = H'_i$.

If $f_i=1$, then let $\widetilde{M}_{i} \subseteq M_i$ be a set of $n^2$
common red neighbours of $\phi(v_{i})$ and $\phi(v_{i+1})$.
We let $H_{i}$ be the subset of $H'_{i}$ formed by the vertices
connected by at least one red edge to a vertex in $\widetilde{M}_{i}$.
If $|H'_{i} \setminus H_{i}| > n^2$, then we have a blue well-split $K_{n^2, n^2}$
with parts $|H'_{i} \setminus H_{i}|$ and $\widetilde{M}_{i}$.
By Corollary~\ref{c:wellSplitBipPW}, this implies a monochromatic non-crossing $G$.
Otherwise we have $|H_{i}| \ge 2mn$.

If $f_i \ge 2$, let $\widehat{M}_{i,1} = M_{i,1} \cap N_{red}(\phi(v_i))$ and
$\widehat{M}_{i,f_i} = M_{i,f_i} \cap N_{red}(\phi(v_{i+1}))$.
For every $j \in [2,f_i-1]$ let $\widehat{M}_{i,j} = M_{i,j}$.
Thus, $|\widehat{M}_{i,j}|\ge 3mn^2$ for every $j \in [f_i]$.

\begin{cl}\lab{cl:starFromTop}
When $f_i \ge 2$, then either there is a monochromatic non-crossing $G$
or there exists a set $H_{i} \subseteq H'_{i}$ of size $2 nm$
such that for every $j \in [f_i]$, every point of $H_{i}$
has $2 m$ red neighbours in $\widehat{M}_{i,j}$.
\end{cl}

\proof
We call a vertex in $H'_{i}$ {\em good}, if it has at least
$2m$ red neighbours in $\widehat{M}_{i,j}$ for every $j \in [f_i]$
and {\em bad} otherwise.
We assume that the number of good vertices is smaller than $2nm$.
The claim will be proven by finding a monochromatic non-crossing copy of $G$.

The number of bad vertices in $H'_{i}$ is at least $|H'_{i}|-2nm\ge 3nm-2nm=nm$.
For each bad vertex $h$, label $h$ by $j$ if $N_{red}(h) \cap \widehat{M}_{i,j} < 2m$.
Since $f_i\le n$, there exists an index $j \in [f_i]$
such that the number of bad vertices labelled $j$ is at least $nm/f_i \ge m \ge n^2$.
Consider the set $\widehat{M}_{i,j}$ and a set
$W$ of $n^2$ bad vertices in $H_{i}$ labelled $j$.
For each $w\in W$, we have $|N_{red}(w) \cap \widehat{M}_{i,j}| \leq 2m$
and we also have $|\widehat{M}_{i,j}| \geq 3mn^2$.
Then by Observation~\ref{o:intersectingSets},
$|\widehat{M}_{i,j} \setminus \bigcup_{w \in W} N_{red}(w)|\ge 3mn^2-2m\cdot n^2 \ge n^2$.
All the edges between $\widehat{M}_{i,j} \setminus \bigcup_{w \in W} N_{red}(w)$ and $W$ are blue.
This implies the occurrence of a well-split blue $K_{n^2,n^2}$
and a monochromatic non-crossing $G$ by Corollary~\ref{c:wellSplitBipPW}.
\qed

Hence we can assume that for every $i$ with $f_i \ge 2$ there is
a set $H_{i}\subseteq H'_{i}$ of $2nm$ points each having $2 m$ red neighbours
in each of $\widehat{M}_{i,1}, \dots, \widehat{M}_{i,f_i}$.
This completes the definition of $H_i$ for every $i \in [\ell-1]$.
Next we take an increasing partition of
each $H_i$ into $|U_i|$ parts, each of size at least $2m$.
By Lemma~\ref{l:longPath}, we either find a blue separable $K_{m,m}$ 
or embed the red path $P_u = (y_i)_{i=1}^{n_2}$ on $S_u$
in such a way that $(\phi(y_i))_{i=1}^{n_2}$ is a red increasing path
with every vertex of $U_i$ mapped on some point of $H_i$.
We consider every star centred at some $u_i$.
If $f_i=1$, then $\phi(u_i)$ has a red neighbour in $\widetilde{M}_{i}$
and we map the only vertex of $Q_i$ on this red neighbour.
If $f_i \ge 2$, let $M'_{i,j}$ be the $2m$ red neighbours of
$\phi(u_i)$ in $\widehat{M}_{i,j}$, for every $j \in [f_i]$.
We assume that there is no blue separable $K_{m,m}$.
Recall that red edges connect all the vertices of $\widehat{M}_{i,1}$ to $\phi(v_i)$ 
and all the vertices of $\widehat{M}_{i,f_i}$ to $\phi(v_{i+1})$.
This fact and Lemma~\ref{l:longPath} imply the existence of a red increasing path 
between $\phi(v_i)$ and $\phi(v_{i+1})$ visiting every ${M}'_{i,j}$ exactly once.
This completes the embedding of the vertices of $G$ on $S_u \cup S_v$
that yields a monochromatic non-crossing graph isomorphic to $G$.
See Fig.~\ref{fig:pw-found}.

\onefigure{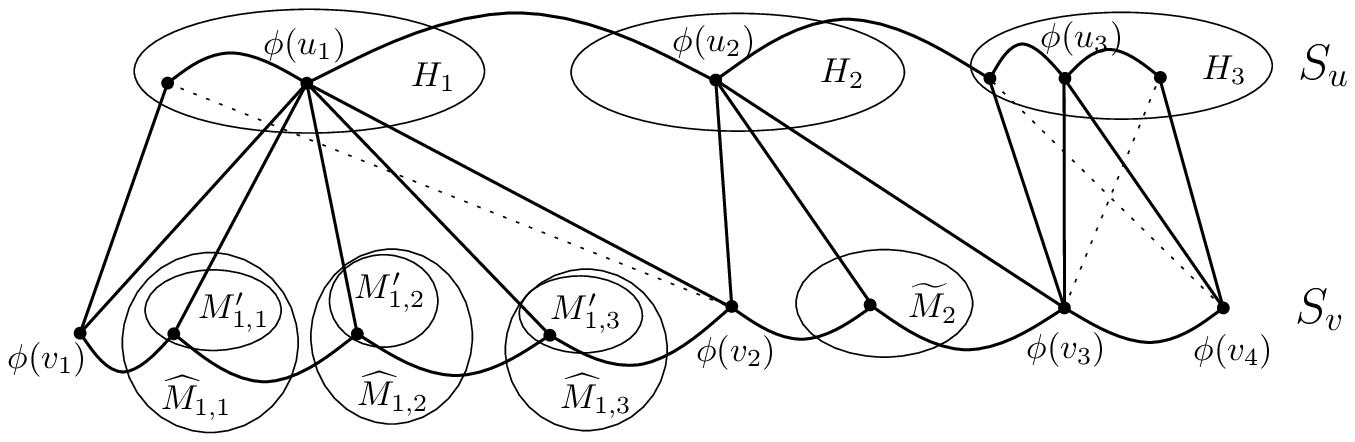}{1}{
Full lines form a red occurrence of the graph from Fig.~\ref{fig:pw-start}.
Dashed lines are other edges known to be red in the colouring $\gamma$. 
}{fig:pw-found}

\section*{Acknowledgements}
The authors would like to thank to Gyula K\'{a}rolyi for introduction to the geometric Ramsey theory
and to Jan Kyn\v{c}l and Martin Balko for discussions about the Ramsey theory of ordered graphs.


\end{document}